# CLUSTERING OF CRITICAL POINTS IN LEFSCHETZ FIBRATIONS AND THE SYMPLECTIC SZPIRO INEQUALITY

V. BRAUNGARDT AND D. KOTSCHICK

ABSTRACT. We prove upper bounds for the number of critical points in semistable symplectic Lefschetz fibrations. We also obtain a new lower bound for the number of nonseparating vanishing cycles in Lefschetz pencils, and reprove the known lower bounds for the commutator lengths of Dehn twists.

## 1. INTRODUCTION

It is a result of Szpiro [16] that in a semistable algebraic family of elliptic curves over $\mathbb{C}P^1$ the number of critical points is bounded above by 6 times the number of singular fibers. In fact, Szpiro considered the arithmetic situation over a number field, and the function field case was just a byproduct of these considerations. Recently, Beauville [3] proved a generalisation of Szpiro's inequality to fibered surfaces, where both the base and the fiber have arbitrary genus.

Amorós et al. [1] gave a group-theoretic proof of Szpiro's inequality for semistable symplectic Lefschetz pencils, and this was extended to certain hyperelliptic semistable symplectic Lefschetz pencils by Bogomolov et al. [4]. Their result is that if all the vanishing cycles are non-separating and the fibration has a topological section, then the number of critical points is bounded above by $4h+2$ times the number of singular fibers, where $h$ is the genus of a smooth fiber.

The purpose of this paper is to prove an inequality of Szpiro type for all semistable symplectic Lefschetz fibrations of arbitrary base and fiber genus, without any assumption on the vanishing cycles and without the hyperelliptic assumption. In the case of pencils our result is that the number of critical points is bounded above by $6h(D-1)$, where $h$ is the genus of a smooth fiber and $D$ is the number of singular fibers. See

Date: January 27, 2003; MSC 2000: 57R17, 57R57, 14H10.

Support from the *Deutsche Forschungsgemeinschaft* is gratefully acknowledged. The authors are members of the *European Differential Geometry Endeavour* (EDGE), Research Training Network HPRN-CT-2000-00101, supported by The European Human Potential Programme.

Theorem 15 and Remark 18 below. In genus one we once more recover the complex function field version of Szpiro's theorem. In higher genus, our inequality is weaker than the one obtained by Beauville [3] in the algebraic situation.

The point of the symplectic Szpiro inequality is that while one can always perturb a Lefschetz fibration in the neighbourhood of a singular fiber so as to make it injective on its critical set, there are global obstructions to the clustering or concentration of critical points in a fiber. These obstructions are essentially the ones discovered by Endo and Kotschick [6] in their proof that the commutator lengths of powers of Dehn twists have linear growth.

In Section 3 we generalise the main result of [6] to relatively minimal Lefschetz fibrations which need not be injective on their critical sets, and then deduce a Szpiro inequality for fibrations over bases of positive genus from this generalisation. We then rederive the lower bounds for the commutator lengths of Dehn twists proved in [6, 8]. There the argument with the Kneser inequality from [9] was combined with the observation that iteration of Dehn twists makes the signatures of Lefschetz fibrations more and more negative. This is clear for separating Dehn twists, but also works for nonseparating Dehn twists by using a handle decomposition, cf. [8]. Here we give a treatment that avoids handle decompositions and derives the growth of the negative definite part of the intersection form from the purely homological Proposition 5 which is standard in algebraic geometry but holds also in our symplectic setup. This makes the argument selfcontained, in particular it is independent of signature calculations for Lefschetz fibrations carried out ad hoc or by using the Meyer signature cocycle.

In Section 4 we generalise results of Li [10] and Stipsicz [14, 15] to relatively minimal symplectic Lefschetz fibrations which need not be injective on their critical sets, and then derive a Szpiro inequality from these generalisations. We also prove that the number $n$ of nonseparating vanishing cycles in a Lefschetz pencil of genus $h$ is no less than $\frac{1}{5}(8h-3)$. The best previously known bound was $n \geq h$, cf. Remark 22. Here too our argument is independent of any signature calculations.

*Acknowledgement:* We are grateful to L. Katzarkov for getting us interested in symplectic analogs of the Szpiro inequality, and for useful conversations. We also like to thank C. Bohr for useful discussions, and A. Stipsicz for reference [15].

## 2. Semistable symplectic Lefschetz fibrations

For definitions and background on differentiable Lefschetz fibrations we refer to [7].

Let $f\colon X \to B$ be an oriented Lefschetz fibration with base genus $g$, fiber genus $h \geq 1$, with $n$ nonseparating and $s$ separating vanishing cycles. We denote by $k = n + s$ the total number of vanishing cycles and by $D$ the number of singular fibers, so that $k \geq D$. We denote by $N$ the total number of components of singular fibers. Note that for $h = 1$ we have $k = N$.

We shall assume throughout that the total space $X$ is a symplectic manifold in such a way that the fibers are symplectic submanifolds. By a theorem of Gompf, cf. [7], this is no restriction if $h \geq 2$.

The Euler characteristic of a Lefschetz fibration is given by

$$\chi(X) = 4(g-1)(h-1) + k \ . \tag{1}$$

The following is elementary:

**Lemma 1.** *If $K$ denotes the canonical class of an almost complex structure associated with the symplectic structure, then*

$$\begin{aligned} K^2 &= 5\chi(X) - 6 + 6b_1(X) - 6b_2^-(X) \\ &= 20(g-1)(h-1) + 5k - 6 + 6b_1(X) - 6b_2^-(X) \ . \end{aligned} \tag{2}$$

**Definition 2.** A symplectic Lefschetz fibration $f\colon X \to B$ is called *semistable* if every 2-sphere component of a singular fiber contains at least two critical points.

The name comes from the fact that if we choose a Riemannian metric on $X$ and consider the induced conformal structures on the fibers, these become semistable algebraic curves if and only if the above topological condition is satisfied.

**Lemma 3.** *A symplectic Lefschetz fibration is semistable if and only if it is relatively minimal.*

*Proof.* Assume it is semistable and suppose a singular fiber has more than one component. Let $S$ be an irreducible component. If $F$ denotes the class of the generic fiber, then $F \cdot S = 0$. Thus

$$S^2 = S \cdot (S - F) \leq -1$$

as the singular fiber is connected, and different components intersect positively unless they are disjoint. If $S \cdot (S - F) = -1$, then by the definition of semistability $S$ is not an embedded sphere. Thus there is no embedded sphere of selfintersection $-1$ in any fiber.

Conversely, if the fibration is relatively minimal, then by the same calculation as above, there can be no spherical component containing only one critical point. □

A variation of this calculation shows the following:

**Lemma 4.** *Let $X$ be a Lefschetz fibration and $F_0$ a singular fiber. Then $H_2(X)$ contains a negative definite subspace for the intersection form spanned freely by the classes of all the components of $F_0$ but one.*

*Proof.* In the case of algebraic surfaces this is well-known as Zariski's lemma; see e. g. [2]. The argument given there goes through in the symplectic situation because the components of singular fibers that we have in the statement intersect each other positively if they are not disjoint. □

**Proposition 5.** *Let $X$ be any Lefschetz fibration. Then*
$$b_2^-(X) \geq 1 + N - D \ . \tag{3}$$

*Proof.* The negative definite subspaces of $H_2(X)$ obtained by applying the preceding lemma to the different singular fibers are mutually orthogonal. Thus the direct sum of all these subspaces is negative definite of dimension $N - D$, and is still orthogonal to the class of a generic fiber, which has zero self-intersection. □

We shall also need the following estimates for the number of components of singular fibers:

**Proposition 6.** *For every Lefschetz fibration we have*
$$N \geq s + D \ , \tag{4}$$

$$N \geq k - (h-1)D \ . \tag{5}$$

*Proof.* The first inequality is immediate from the definition.

To prove (5), let $\Sigma$ be the union of the singular fibers. We can calculate the Euler characteristic of $\Sigma$ by comparing the singular fibers to a generic fiber $F$
$$\chi(\Sigma) = D\chi(F) + k = -2(h-1)D + k \ ,$$
or by summing over all components $C_1, \ldots, C_N$ of the singular fibres
$$\chi(\Sigma) = \sum_{i=1}^N \chi(C_i) - k \leq 2N - k \ .$$
Thus we have $N \geq k - (h-1)D$. □

*Remark* 7. For fibrations with only stable fibers we have $k \leq 3(h-1)D$ and $N \leq 2(h-1)D$ for purely topological reasons.

## 3. Lefschetz fibrations over bases of positive genus

In this section we prove the main technical result about relatively minimal Lefschetz fibrations and derive inequalities of Szpiro type under the assumption that the base has positive genus.

**Theorem 8.** *Let $X$ be a connected smooth closed oriented 4-manifold and $f\colon X \to B$ a relatively minimal symplectic Lefschetz fibration with fiber genus $h \geq 1$ and base genus $g \geq 1$ having $k$ vanishing cycles, $D$ singular fibers and $N$ irreducible components of singular fibers. Then*

$$5k + 6(3h-1)(g-1) \geq 6(N-D) \ . \tag{6}$$

*Proof.* As $X$ is assumed to be relatively minimal, the positivity of the base genus implies that $X$ is minimal and not ruled, because any pseudo-holomorphic sphere in $X$ would have to be contained in a fiber. Thus Liu's extension [11] of Taubes's results [17] implies $K^2 \geq 0$, which we can write as

$$b_2^+(X) \geq \frac{1}{5}(b_2^-(X) + 4b_1(X) - 4) \ . \tag{7}$$

Using (3) and $b_1(X) \geq 2g \geq 2$, we obtain

$$b_2^+(X) \geq 1 + \frac{1}{5}(N-D) \ .$$

As the claim (6) is trivial for $N = D$, we may assume $N - D \geq 1$, and therefore $b_2^+(X) \geq 2$.

As $X$ is minimal with $b_2^+(X) \geq 2$, we can use the result of Taubes [17] to obtain a symplectically embedded surface $\Sigma \subset X$ representing the canonical class $K$ of $X$. It may be disconnected, but because $X$ is minimal, $\Sigma$ has no spherical component. In the argument below we will tacitly assume that it is connected. In the general case the same argument works by summing over the components.

The genus of $\Sigma$ is given by the adjunction formula $g(\Sigma) - 1 = K^2$. The fibration $f$ induces a smooth map $\Sigma \to B$ of degree $d$ equal to the algebraic intersection number of $\Sigma$ with a fiber. This is calculated from the adjunction formula applied to a smooth fiber $F$, which is a symplectic submanifold:

$$d = \Sigma \cdot F = K \cdot F = 2h - 2 \ . \tag{8}$$

Thus Kneser's inequality $g(\Sigma) - 1 \geq |d|(g(B) - 1)$ gives:

$$K^2 \geq 2(h-1)(g-1) \ .$$

Combining this with (2) and estimating $K^2$ from above using $b_2^-(X) \geq 1 + N - D$ and $b_1(X) \leq 2g + 2h$ we obtain

$$6(N - D) \leq 18(h - 1)(g - 1) - 12 + 12g + 12h + 5k \ .$$

Pulling the fibration back to large degree covers of the base $B$ and applying the above inequality, we finally obtain (6). □

**Corollary 9.** *Let $X$ be a connected smooth closed oriented 4-manifold and $f \colon X \to B$ a relatively minimal symplectic Lefschetz fibration with fiber genus $h \geq 1$ and base genus $g \geq 1$ having $s$ separating and $n$ non-separating vanishing cycles, $D$ singular fibers and $N$ irreducible components of singular fibers. Then the following inequalities hold:*

$$(9) \qquad s \leq 6(3h - 1)(g - 1) + 5n \ ,$$

$$(10) \qquad k \leq 6(3h - 1)(g - 1) + 6hD \ ,$$

$$(11) \qquad N \leq 6(3h - 1)(g - 1) + (5h + 1)D \ .$$

*Proof.* The first claim follows from (6) using $N \geq s + D$. The second and third claim follow similarly using $k \leq N + (h - 1)D$. □

*Remark* 10. The inequality (9) was originally proved by Endo and Kotschick [6] under the assumption $k = D$.

*Remark* 11. Note that for $h = 1$ we have $s = 0$ and $k = N$. From (10) or (11) we obtain

$$k = N \leq 12(g - 1) + 6D \ ,$$

which is a generalisation of the Szpiro inequality to semistable symplectic Lefschetz fibrations over bases of positive genus. For $h \geq 2$, either (10) or (11) can be regarded as a Szpiro-type inequality.

We now apply the previous discussion to give a new proof for the known lower bounds for the commutator lengths of powers of Dehn twists in mapping class groups.

**Theorem 12.** *Let $a$ be a homotopically nontrivial simple closed curve on a surface $F$ of genus $h \geq 2$, and let $t_a$ be the corresponding Dehn twist. Suppose that $t_a^k$ with $k > 0$ can be written as a product of $g$ commutators. Then*

$$(12) \qquad g \geq 1 + \frac{k}{6(3h - 1)} \ .$$

*Proof.* Consider the standard holomorphic Lefschetz fibration over the 2-disk $D^2$ with precisely one singular fiber $F_0$ with vanishing cycle $a$. Pulling back under the base change $z \mapsto z^k$ and taking the minimal resolution, we obtain a holomorphic Lefschetz fibration over $D^2$ with only one singular fiber having $k$ vanishing cycles which are parallel copies of $a$. The monodromy of this fibration around the boundary of the disk is $t_a^k$. If this can be expressed as a product of $g$ commutators, then we can find a smooth surface bundle with fiber $F$ over a surface of genus $g$ with one boundary component and the same restriction to the boundary. Let $X$ be the Lefschetz fibration over the closed surface $B$ of genus $g$ obtained by gluing together the two fibrations along their common boundary.

By construction, $X$ is symplectic and relatively minimal, so that we can apply Theorem 8. We have $D = 1$ and $N = k + 1$ if $a$ is separating and $N = k$ if $a$ is nonseparating. Theorem 8 gives $k \leq 6(3h - 1)(g - 1) + c$ with $c = 0$ or $c = 6$ depending on whether $a$ is separating or not. In the latter case by pulling back the fibration to large degree coverings of the base we also obtain $k \leq 6(3h - 1)(g - 1)$ as claimed. □

*Remark* 13. The inequality (12) was originally proved by Endo and Kotschick [6] under the assumption that $a$ is separating. In the context of Lefschetz fibrations with $k = D$ they used the separating assumption to conclude $b_2^-(X) \geq k + 1$ whenever $a$ occurs $k$ times as a vanishing cycle, as every separating vanishing cycle makes a negative contribution to the signature. Korkmaz [8] then observed that using a handle decomposition one still gets $b_2^-(X) \geq k$ in the nonseparating case, so that the argument goes through. Phrased as above, the proof works directly for both cases, as the required lower bound for $b_2^-(X)$ arises from the homological argument in Proposition 5.

## 4. Lefschetz pencils

In this section we consider the case of pencils, i. e. Lefschetz fibrations over the 2-sphere. We shall assume that our pencils are nontrivial, meaning that they have at least one critical point each.

As in [10, 15], the case of ruled surfaces has to be considered separately.

**Proposition 14.** *Suppose $X$ is the blowup in $b$ points of a 2-sphere bundle over a surface of genus $a$, and that $X$ admits a nontrivial relatively minimal symplectic Lefschetz pencil with fiber genus $h \geq 1$ having $s$ separating and $n$ nonseparating vanishing cycles, $D$ singular fibers* *and $N$ irreducible components of singular fibers. Then the following inequalities hold:*

$$(13) \qquad k \geq 2h - 2 + \frac{3}{2}(N - D) ,$$

$$(14) \qquad n \geq 2h - 2 + \frac{1}{2}(N - D) .$$

*Proof.* By the assumptions on $X$, we have $b_2^+ = 1$, $b_2^- = 1 + b$ and $b_1 = 2a$.

Stipsicz [14] proved that for a nontrivial relatively minimal Lefschetz pencil $K^2 \geq 4(1-h)$. His proof was written under the assumption that Lefschetz pencils are injective on their critical sets, but this can be achieved by perturbation, and the inequality involves only topological invariants which do not change under perturbation (unlike $D$ and $N$), so the inequality is true in our case. Substituting the above numbers into it, we obtain

$$(15) \qquad 4a \leq 2 + 2h - \frac{1}{2}b .$$

Computing the Euler characteristic of $X$ in two different ways we see $4(1-a) + b = \chi(X) = 4(1-h) + k$, and so

$$k = b + 4h - 4a \geq 2h - 2 + \frac{3}{2}b = 2h - 2 + \frac{3}{2}(b_2^- - 1) \geq 2h - 2 + \frac{3}{2}(N - D) ,$$

where we have used first (15) and then Proposition 5. Thus we have proved the first claim. The second one follows from the first using Proposition 6. □

In general we have:

**Theorem 15.** *Let $X$ be a connected smooth closed oriented 4-manifold and $f\colon X \to S^2$ a nontrivial relatively minimal symplectic Lefschetz pencil with fiber genus $h \geq 1$ having $s$ separating and $n$ nonseparating vanishing cycles, $D$ singular fibers and $N$ irreducible components of singular fibers. Then the following inequalities hold:*

$$(16) \qquad 5k \geq 6h + 6(N - D) ,$$

$$(17) \qquad 5n \geq 6h + s ,$$

$$(18) \qquad k \leq 6h(D - 1) ,$$

$$(19) \qquad N \leq (5h+1)(D-1) - (h-1) \ .$$

*Proof.* First we fiber sum $X$ with a genus $h$ bundle over the 2-torus, and apply (9) to the resulting Lefschetz fibration. This shows that, as $X$ is nontrivial, we must have $n > 0$. Therefore $b_1(X) \leq 2h - 1$.

Assuming first that $X$ is not rational or ruled, we would like to use Taubes's result [17] as extended by Liu [11] to obtain $K^2 \geq 0$. The problem with this argument is that in the case of base genus zero, relative minimality does not imply minimality, so that Taubes's result is not available. However, if $X$ is not rational or ruled, we have Li's inequality [10]

$$(20) \qquad K^2 \geq 2 - 2h \ .$$

Li's argument assumes that the Lefschetz pencil is injective on its critical set, but this can be achieved by perturbation without affecting the inequality.

Using this, we proceed as in the proof of Theorem 8. Combining (20) with (2) we obtain (16) by using $g = 0$, $b_1(X) \leq 2h - 1$ and $b_2^-(X) \geq 1 + N - D$.

It remains to deal with the case that $X$ is a (blowup of a) ruled surface. If $h \geq 3$, then (16) follows from (13).

If $h = 2$, suppose we have a fibration satisfying (13) but failing (16). Then we conclude $k = 2$, giving $\chi(X) = 4(1-h) + k = -2$. On the other hand, in the notation of Proposition 14 we have $a \leq \frac{1}{2}h + \frac{1}{2} - \frac{1}{8}b$, which gives $a \leq 1$. We conclude $\chi(X) = 4(1-a) + b \geq 0$, which is a contradiction.

If $h = 1$ for a fibration $X$ satisfying (13) but failing (16), then we conclude $k \leq 3$ and therefore $\chi(X) = k \leq 3$. As above we also obtain $a \leq 1$. If $a = 0$, then $\chi(X) = 4$, which is a contradiction. If $a = 1$, then (15) gives $b = 0$ and therefore $k = \chi(X) = 0$. But then the fibration is trivial.

Thus (16) is proved in all cases. Once we have (16), the other inequalities follow as in the proof of Corollary 9. $\square$

*Remark* 16. For non-ruled total spaces, the inequality (17) was proved by Li [10] under the assumption $k = D$.

*Remark* 17. If $k \neq 0$, we conclude $D \geq 2$ from (18). Thus, the critical points of a non-trivial Lefschetz pencil can never be concentrated in a single fiber, compare [12].

*Remark* 18. Note that for $h = 1$ we have $k = N$ and both (18) and (19) reprove Szpiro's inequality [16] and give extensions to pencils with $h >$

1. Concerning (18), note that Bogomolov et al. [4] proved $k \leq (4h + 2)D$ for hyperelliptic fibrations under the additional assumption that all vanishing cycles are nonseparating and that the fibration admits a topological section.

*Remark* 19. In the proof of Theorem 15 we have used $b_1(X) \leq 2h - 1$. If we actually know the first Betti number, then we obtain better inequalities.

The following theorem gives new bounds on the number of nonseparating vanishing cycles in Lefschetz pencils.

**Theorem 20.** *Let $X$ be a connected smooth closed oriented 4-manifold and $f\colon X \to S^2$ a nontrivial relatively minimal symplectic Lefschetz pencil with fiber genus $h \geq 1$ having $s$ separating and $n$ nonseparating vanishing cycles, $D$ singular fibers and $N$ irreducible components of singular fibers. If $X$ is not rational or ruled, then the following inequalities hold:*

$$5k \geq 8h - 3 + 5(N - D) , \tag{21}$$

$$5n \geq 8h - 3 . \tag{22}$$

*Proof.* The assumption that $f$ is relatively minimal again does not imply that $X$ is minimal, in which case we would have $K^2 \geq 0$. However, $K^2 \geq 0$ can be rewritten as

$$5b_2^+ - 4b_1 + 4 \geq b_2^- . \tag{23}$$

Blowing up or down does not change the left hand side of this inequality, so, regardless of whether $X$ is minimal or not, the fact that it is not rational or ruled implies

$$5b_2^+(X) - 4b_1(X) + 4 \geq 0 \tag{24}$$

because of [17, 11]. Using (1) and

$$b_2^+(X) = b_2(X) - b_2^-(X) = \chi(X) - 2 + 2b_1(X) - b_2^-(X) ,$$

we obtain

$$5k \geq 20h - 14 - 6b_1(X) + 5b_2^- . \tag{25}$$

As in the proof of Theorem 15 we have $b_1(X) \leq 2h - 1$. From Proposition 5 we have $b_2^-(X) \geq N - D + 1$. Thus (25) gives (21).

Using $N - D \geq s$, we obtain (22) from (21). □

Finally we extend (22) to all Lefschetz pencils:

**Theorem 21.** *If $f\colon X \to S^2$ is a nontrivial Lefschetz pencil of genus $h \geq 1$, then $f$ has at least $\frac{1}{5}(8h-3)$ nonseparating vanishing cycles.*

*Proof.* Clearly it suffices to prove the relatively minimal case. If $X$ is not rational or ruled, then the result is part of the previous Theorem, cf. (22). If $X$ is rational or ruled, we have (14), which is enough as long as $h \geq 4$. On the other hand, if $h \leq 2$, then (17) implies the claim. Thus it remains to deal with the case $h = 3$ for ruled manifolds.

If $h = 3$ and $X$ is the blowup of a ruled surface satisfying (14) but failing $n \geq \frac{1}{5}(8h-3)$, then again $s \leq N - D = 0$, and $k = n = 4$. The same arguments as in the proof of Theorem 15 above then give $b = 0$ and $a = 2$. Thus we have a ruled surface without any blowups over a genus 2 surface, in which the generic fiber $F$ of the Lefschetz fibration represents a homology class of zero selfintersection. As $F$ is a symplectic submanifold, the proof of the Thom conjecture implies that $F$ has minimal genus in its homology class. We now show that this leads to a contradiction.

If $X$ is the product ruled surface $S^2 \times \Sigma_2$, then the homology classes of zero selfintersection are the multiples of the two factors. The multiples of the first factor are all represented by spheres, so $F$ could only be a multiple of the second factor. As it has larger genus than the second factor, it would represent $d$ times the second factor with $|d| \geq 2$. But the projection to the second factor induces a map of degree $d$ from $F$ and from the normalisation of a singular fiber to $\Sigma_2$. As the normalisation of a singular fiber has genus $\leq 2$, it does not map to $\Sigma_2$ with degree $\geq 2$.

For the nontrivial ruled surface, $F$ can also not be a multiple of the class of the 2-sphere in the ruling, and must therefore be a class which has intersection number $d \neq 0, \pm 1$ with the class of the sphere in the ruling. Thus it maps to the base surface of genus 2 with degree $\geq 2$, and we obtain the same contradiction as in the product case. □

*Remark* 22. Theorem 6.1 of Stipsicz's paper [15] gives the lower bound $\frac{1}{5}(8h-4)$ for the total number of critical points in nontrivial Lefschetz pencils (which are injective on their critical sets). However, the proof given there does not seem to cover the case of ruled surfaces with fibrations of small fiber genus, because it appeals to Theorem 1.1 which gives nothing for fiber genus $\leq 5$, for example. More importantly, Stipsicz's argument appeals to Theorem 1.4 of that paper, whose proof is incomplete, because it uses a statement about the fundamental group of a Lefschetz fibration which is not known to be true. In fact, it is false if one considers achiral Lefschetz fibrations, compare [1], p. 503.

In any case, the important difference between the result of Stipsicz and ours is the fact that we obtain a bound growing with $\frac{8}{5}h$ for the number of nonseparating vanishing cycles, rather than the total number of separating and nonseparating ones. In the nonseparating case the best previous result is that of Li [10], which is $n \geq h$.

*Remark* 23. Combining Theorems 15 and 21, the number of nonseparating vanishing cycles in a Lefschetz pencil is bounded below by

$$5n \geq ts + (8 - 2t)h + 3(t - 1) \tag{26}$$

for all $t \in [0, 1]$.

## 5. Final comments

As in [9, 6], the arguments of this paper rely on the work of Taubes [17] in Seiberg-Witten theory, showing that for a minimal symplectic 4-manifold with $b_2^+ > 1$, the canonical class is represented by a symplectically embedded surface without spherical components. Recently, the methods pioneered in Donaldson's work on Lefschetz pencils have been applied successfully to reprove Taubes's result. Donaldson and Smith [5] did this under the additional assumption $b_2^+ - b_1 > 1$. This is not sufficient for our purposes, but, as explained by Smith at the end of [13], the arguments of [5] can be pushed to cover all cases where $b_2^+ > 2$.

This means that as long as we work with manifolds with $b_2^+ > 2$, our results can be proved independently of gauge theory. In Section 3 we appealed to Liu's work [11] in gauge theory for the case $b_2^+ = 1$, but this can be avoided. Therefore, our results on fibrations over bases of positive genus do not require any input from gauge theory. In the case where the base genus is zero, it is possible that $b_2^+ = 1$ or 2, but if we exclude those cases, then the proofs can be based on [5, 13] instead of [17, 11].

MATHEMATISCHES INSTITUT, LUDWIG-MAXIMILIANS-UNIVERSITÄT MÜNCHEN, THERESIENSTR. 39, 80333 MÜNCHEN, GERMANY

*E-mail address*: `Volker.Braungardt@mathematik.uni-muenchen.de`

MATHEMATISCHES INSTITUT, LUDWIG-MAXIMILIANS-UNIVERSITÄT MÜNCHEN, THERESIENSTR. 39, 80333 MÜNCHEN, GERMANY

*E-mail address*: `dieter@member.ams.org`